\documentclass[12pt,a4paper]{amsart}
\usepackage{amsmath,amsthm,amscd,amssymb}
\usepackage{latexsym}
\usepackage{epic, eepic}
\def\squarebox#1{\hbox to #1{\hfill\vbox to #1{\vfill}}}

\newcommand{\cz}{{\mathbb C}}
\newcommand{\nz}{{\mathbb N}}
\newcommand{\qz}{{\mathbb Q}}
\newcommand{\rz}{{\mathbb R}}
\newcommand{\zz}{{\mathbb Z}}

\def\qed{\hbox {\hskip 1pt \vrule width 4pt height 6pt depth 1.5pt
        \hskip 1pt}}% Beweisende
\def\canform{canonical transformation{}}

\def\fourior{Fourier integral operator{}}

\def\neigh{neighborhood{}}

\def\pseudor{pseudodifferential operator{}}

\def\wrt{with respect to}

\def\Re{{\rm Re\,}}
\def\Im{{\rm Im\,}}
\newtheorem{theorem}{Theorem}

\newtheorem{lemma}{Lemma}

\title[Trapping point resonances]
{An Inverse problem for trapping point resonances.}
\author{Alexei Iantchenko}
\address{Institute of Mathematics and Physics\\
University of Wales, Aberystwyth\\
Penglais, Ceredigion UK SY23 3BZ}  \email{ aii@aber.ac.uk}
\address{On leave from Malm\"o H{\"o}gskola,
Sweden}
\date{\today}
\begin{document}

\bibliographystyle{plain}



\begin{abstract}
We consider  semi-classical Schr{\"o}dinger operator $
P(h)=-h^2\Delta +V(x)$ in ${\mathbb R}^n$ such that the analytic
potential $V$  has a non-degenerate critical point $x_0=0$  with
critical value $E_0$ and we can define resonances in  some fixed
neighborhood of $E_0$ when $h>0$ is small enough.  If the
eigenvalues of the Hessian are $\zz$-independent the resonances in
$h^\delta$-neighborhood of $E_0$ ($\delta >0$) can be calculated
explicitly as the eigenvalues of the semi-classical Birkhoff normal
form.

Assuming that potential is symmetric with respect to reflections
about the coordinate axes  we show that the classical Birkhoff
normal form  determines the Taylor series of the potential at $x_0.$
As a consequence,    the resonances in a $h^\delta$-neighborhood of
$E_0$ determine the first $N$ terms in the Taylor series of $V$ at
$x_0.$

The proof uses the recent  inverse spectral results of V. Guillemin
and A. Uribe.

\end{abstract}

\pagestyle{myheadings} \maketitle

\section{Introduction.}
%In the present note we show that the method of V. Guillemin and A.
%Uribe in \cite{GuilleminUribe2007} on spectral semi-classical
%inverse problem can be extended to the resonances.

We consider the semi-classical Schr{\"o}dinger operator
\begin{equation}\label{1.1}
P=P(h)=-h^2\Delta +V(x),\,\,x\in\rz^n,\end{equation}   with the
symbol $p(x,\xi)=\xi^2+V(x).$

If the spectrum of (\ref{1.1}) is discrete near some energy $E$ and
real-valued potential  $V$ is smooth then it is known
(\cite{Guillemin1996}, \cite{Zelditch2004},
\cite{IantchenkoSjostrandZworski2002}) that the spectrum of $P(h)$
in a small fixed \neigh{} of $E$ as $h\rightarrow 0$ determines the
Birkhoff normal form of the Hamiltonian $p(x,\xi)=\xi^2 +V(x).$ In
\cite{GuilleminUribe2007} it was shown that the classical Birkhoff
normal form of $p(x,\xi)$ at a non-degenerate minimum $x_0$ of $V$
determines the Taylor series of the potential provided the
eigenvalues of the Hessian are linearly independent over $\qz$ and
$V$ satisfies a symmetry condition near $x_0.$ This result was
applied  to prove that the low-lying eigenvalues of the
semi-classical operator $P$ determine the Taylor series of the
potential at $x_0.$ In this note we study the similar question for
the resonances. In \cite{Zworski2007} it was indicated how the
inverse spectral results based on wave invariants translates to
inverse results for resonances (see also \cite{Zelditch2002}).  We
consider a special situation as in \cite{Sjostrand1986} and
\cite{KaidiKerdelhue2000}, when the resonances can be calculated
explicitly as the eigenvalues of the semi-classical Birkhoff normal
form.

%In the several papers ...it was shown that the Birkhoff normal form
%can be recovered from semi-classical spectral invariants.

%describes behaviour of the classical dynamics near the point of
%stable equilibrium. Its quantization was successfully applied to the
%calculation of the semi-excited states for the Schr{\"o}dinger
%equation with the potential whose Taylor series at $0$ starts with
%quadratic terms. Later several generalizations of the semi-classical
%or quantum Birkhoff normal form (qBnf)  have appeared and applied to
%the semi-classical inverse problems. Namely it was shown that the
%qBnf is spectral invariants. Using the Helffer-Sj{\"o}strand theory
%of resonances the qBnf could also be used to describe the resonances
%near some fixed energy Kaidi-Kerdelhue, Sj, and scattering resonances for
%Helmholtz equation outside two strictly convex obstacles generated
%by the trapped billiard trajectory between the obstacles. In the
%present note we will show that the recent detailed inverse
%semi-classical spectral results can be generalized to the
%resonances.

 We suppose  that general assumptions of
Helffer-Sj{\"o}strand in \cite{HelfferSjostrand1986} are fulfilled
so that we can define resonances in some fixed \neigh{} of $E_0\in\rz$ when $h>0$ is small enough.

We suppose also that $V$ is  analytic potential, which extends to a
holomorphic function in a set
$$\{ x\in\cz^n;\,\,|\Im x| <\frac{1}{C}\langle\Re x\rangle\}$$ with
$V(x)\rightarrow 0,$ when $x\rightarrow \infty$ in that set. Here
$\langle s\rangle=(1+|s|^2)^{1/2}.$

%The resonances of $P$ can be defined in an angle
%$\{z\in\cz;\,\,-2\theta_0 <\arg\, z \leq 0\}$ for some fixed
%$\theta_0
%>0$ as the  eigenvalues of
%$P{_{|e^{i\theta_0}\rz^n}}.$

%Using the complex deformation argument the resonances can be defined
%as the poles of the meromorphic continuation of the resolvent
%$$R(z) =(P-z)^{-1}:\,\,L_{\rm comp}^2(\rz^{n})\mapsto L_{\rm loc}^2(\rz^{n})$$ from
%$\Im z>0$ across $]0,\infty[$ to $e^{-i[0,2\theta_0[}]0,\infty[.$
%The poles are the resonances and the multiplicity of a resonance is
%given by $$m(z_0)=-{\rm rang}\int_{\gamma}R(z)dz,$$ where $\gamma$
%is a small circle centered at $z_0.$

 We will use notation ${\rm neigh}\,({E,\rz})$ or  ${\rm neigh}(E)$ for a real \neigh{} of a $E\in\rz.$

Following \cite{GerardSjostrand1987} the trapped set $K(E_0)$   is
$K
(E_0)=\{\rho\in p^{-1}(E_0);\,\,\exp tH_p(\rho)\not\rightarrow
\infty,\,\,t\rightarrow\pm\infty\},$ which is the union of trapped
trajectories in $p^{-1}(E_0).$ Here $H_p$ is the Hamilton field of
$p(x,\xi).$

We assume that the union of trapped trajectories in $p^{-1}(E_0)$ is
just the point $(0,0):$
\begin{equation}\label{trapped point}
K(E_0)=(0,0).
\end{equation}
Then $0$ is a unique critical point of $V$  with critical value
$E_0.$ We suppose that $0$ is non-degenerate critical point of $V$
with signature $(n-d,d):$
 $$V(0)=E_0,\,\,V'(0)=0,\,\,{\rm sgn}\,V''(0)=(n-d,d),$$ so that
$V''(0)$ is non-degenerate and \begin{equation}\label{V}
V(x)=E_0+\sum_{j=1}^{n-d} u_j^2x_j^2 -\sum_{j=n-d+1}^n u_j^2x_j^2
+{\mathcal O}(|x|^3).\end{equation}

Kaidi and Kerdelhue showed in \cite{KaidiKerdelhue2000} how to adapt the Helffer-Sj{\"o}strand theory
and realize $P=-h^2\Delta +V(x)$ as acting in $H(\Lambda)$-spaces, where
$\Lambda\subset\cz^{2n}$ is an IR-manifold which coincides with
$T^*(\rz^{n-d}\oplus e^{i\pi/4}\rz^{d})$ near $(0,0)$ and has the
property that $\forall\epsilon >0,\,\,\exists\delta >0$ such that
$(x,\xi)\in\Lambda,\,\,{\rm dist}\,((x,\xi),(0,0)) >\epsilon$
$\Rightarrow$ $|p(x,\xi) -E_0| >\delta.$

Then resonances can essentially (modulo an argument using a Grushin
reduction) be viewed as an eigenvalue problem for $P$ after the
complex scaling $ x_j=e^{i\pi/4} \tilde{x}_j,$ $\tilde{x}_j\in\rz,$
$n-d+1\leq j\leq n.$

We suppose also that the coefficients $u_j$ in (\ref{V}) satisfy
non-resonance condition:
\begin{equation}\label{non-res}
\sum_{j=1}^n k_j u_j=0,\,\,
k_j\in\zz\,\,\Rightarrow\,\,k_1=k_2=\ldots=k_n= 0.
\end{equation}

Under these assumptions a result of Kaidi and Kerdelhue
\cite{KaidiKerdelhue2000}  gives all resonances in a disc $D(E_0,
h^\delta)$ of center $E_0$ and radius $h^\delta.$ Here $\delta$ can
be any fixed constant and $h>0$ is small enough depending on
$\delta.$ We cite this result later in Theorem \ref{KaKe}.

The consequence of the main result of this note is the following:
\begin{theorem}\label{TheoremMain}
Assume $V$ is symmetric \wrt{} reflections about the coordinate
axes, i.e. for any choice of signs
\begin{equation}\label{symmetric}
V(x_1,\ldots,x_n)=V(\pm x_1,\ldots,\pm x_n).\end{equation} In
addition, assume that
\begin{equation}\label{expansion}V(x)=E_0+\sum_{j=1}^{n-d} u_j^2x_j^2
-\sum_{j=n-d+1}^n u_j^2x_j^2 +{\mathcal O}(|x|^4),\end{equation}
where  $u_1,\ldots,u_n$ are the positive numbers satisfying
(\ref{non-res}).
%and two families $\{ u_j\}_{j=1}^{n-d}$  and $\{ u_j\}_{j=n-d+1}^{n}$ are $\zz$-independent.

%\begin{equation}\label{non-resonance}
%\sum_{j=1}^{n-d}\alpha_j u_j=0,\,\,\sum_{j=n-d+1}^{n}\alpha_j
%u_j=0,\mbox{for
%some}\,\,\alpha\in\zz^n\,\,\Rightarrow\,\,\alpha\equiv 0.
%\end{equation}

Then, given $N>0$ there exists a $\delta >0$ such that the
resonances in $D(E_0, h^\delta)$ for $0<h<h_0,$ determine the first
$N$ terms in the Taylor series of $V$ at zero.

\end{theorem}
In dimension $n=1,$ $d=1,$ resonances generated by the maximum of
the potential (barrier top resonances)  are of the form $ \simeq
V(0)-ih(-V''(x_0)/2)^{1/2} (2k+1)+\ldots,$ $k=0,1,\ldots,$
$V(0)=E_0.$
 Yves Colin de Verdi{\`e}re and Victor Guillemin have recently
shown in \cite{CdVGuillemin2008} that one can drop the condition
that the potential is even. Namely instead of (\ref{symmetric}) and
(\ref{expansion}) it is enough to suppose that in the expansion
$V(x)=E_0-ux^2+\sum_{j=3}^\infty a_jx^j$ the coefficients $u>0$ and
$a_3$ do not vanish. Then all $a_j$'s  are determined from the
coefficients of the quantum Birkhoff normal form once we have chosen
the sign of $a_3.$ The classical Birkhoff normal form along is not
enough to recover the potential.

In dimension $n=2,$ $d=1,$  Sj{\"o}strand (see \cite{Sjostrand2003})
showed that the saddle-point resonances are given by the eigenvalues
of the Birkhoff normal form in the whole $h$-independent \neigh{} of
$E_0.$ Thus the full Taylor series of $V$ is determined and using
the analyticity, the full potential can be recovered from the
resonances.

To prove Theorem \ref{TheoremMain} we use that under non-resonance
condition (\ref{non-res}) the Schr{\"o}dinger operator $P$ can be
transformed in the semi-classical or quantum Birkhoff normal form
(see \cite{Sjostrand1992})
\begin{equation}\label{Bnf} E_0+\tilde{P}\equiv U^*PU,
\end{equation}
where $U$ is analytic unitary \fourior{} microlocally defined near
$(0,0)$ and  $\tilde{P}$ is \pseudor{} with the symbol
\begin{equation}\label{F} F\sim\sum_{j=0}^\infty
h^jF_j(\imath_1,\ldots,\imath_{n-d} , \jmath_{n-d+1},\ldots,
\jmath_{n} ), \,\,\imath_j=\xi_j^2
+x_j^2,\,\,\jmath_j=\xi_{j}^2-x_{j}^2,
\end{equation}
 with $F_j$ analytic and principal symbol
\begin{equation}\label{F0}
F_0=\sum_{j=1}^{n-d} u_j\imath_j +\sum_{j=n-d+1}^n u_j\jmath_j
+{\mathcal O}(|(\imath,\jmath)|^2).
\end{equation}
The equivalence relation $\equiv$ means  to infinite order at
$(0,0)$ (see \cite{IantchenkoSjostrand2002}).

The result of \cite{KaidiKerdelhue2000} shows
 that, modulo error
terms of order ${\mathcal O}(h^\infty),$ the resonances of $P$ in
$h^\delta$ \neigh{} of $E_0$ are approximated by  the eigenvalues of
its quantum Birkhoff normal form at $(0,0)$ after  the complex
scaling $ x_j=e^{i\pi/4} \tilde{x}_j,$ $\tilde{x}_j\in\rz,$
$n-d+1\leq j\leq n,$ namely $$\tilde{F} \sim \sum_{j=0}^\infty
h^jF_j(\imath_1,\ldots,\imath_{n-d},\frac{1}{i}\tilde{\imath}_{n-d+1},\ldots,\frac{1}{i}\tilde{\imath}_{n}),$$
where $F$  is as in (\ref{F}) and
$\frac{1}{i}\tilde{\imath}=\frac{1}{i}(\tilde{\xi}_j^2+\tilde{x}_j^2)=\xi_j^2-x_j^2,$
$\xi_j=e^{-i\pi/4}\tilde{\xi}_j,$ $x_j=e^{i\pi/4}\tilde{x}_j.$ We denote
$\tilde{F}_j(\imath_1,\ldots,\imath_{n-d},\tilde{\imath}_{n-d+1},\ldots,\tilde{\imath}_{n})=
F_j(\imath_1,\ldots,\imath_{n-d},\frac{1}{i}\tilde{\imath}_{n-d+1},
\ldots,\frac{1}{i}\tilde{\imath}_{n}).$

%The proof is given in Section \ref{Section1} by applying the method
%of GuUr \cite{GuilleminUribe2007} to the complex scaled operator.

\begin{theorem}[Kaidi-Kerdelhue]\label{KaKe} The resonances of $P$ in rectangle $]E_0-\epsilon_0,E_0+\epsilon_0[
-i[0,h^{\delta}]$ are simple labeled by $k\in\nz^n$ and of the form
$$E_0+\sum_{j=0}^\infty h^j\tilde{F}_j((2k_1+1)h,\ldots,(2k_n+1)h)$$ where
$$\tilde{F}_j\in C^\infty({\rm neigh}(0)),\,\, \tilde{F}_0(\imath)=\sum_{j=1}^{n-d}
u_j\imath_j -\sum_{j=n-d+1}^n iu_j\imath_j
+\mathcal{O}(|\imath |^2),\,\,\tilde{F}_1(\imath)=V(0)-E_0=0.$$
\end{theorem}

The main result of this note is the following:
\begin{lemma}\label{lemma1} Assume  (\ref{non-res}), (\ref{symmetric}) and
(\ref{expansion}).
%, where the positive numbers $u_1,\ldots,u_2$ are
%linearly independent over the rationals.
Then the classical Birkhoff normal form $F_0$ determines the
Taylor series of $V$ at the origin.

\end{lemma}
We show in Section \ref{Section1} how this lemma follows from
\cite{GuilleminUribe2007}. The main idea of the proof is that the
complex scaling reduces the principle symbol of $P$ to the form
$H(x,\xi)=\sum_{j=1}^n\omega_j(\xi_j^2+x_j^2)+{\mathcal O}(|x|^3)$
which is similar to the Hamiltonian considered in
\cite{GuilleminUribe2007} with the only difference that coefficients
$\omega_j$ for $n-d+1\leq j\leq n$ are complex numbers. We show in
Section \ref{Section1} that the method of Guillemin and Uribe can
still be applied.\\ \\

%Using decomposition of the Bnf as in \cite{Sjostrand2003} we get
%resonances of the form
%$$E_0+K(0;h)+G((2k+1)\mu)+\sum_{j=1}^\infty h^{j+1}
%(k_j(2k+1)+h_{j+1}(2k+1)),$$ where $k_j$ is polynomial of degree
%$\leq j,$ $k_j(0)=0,$ $h_j(y)=\sum_{|\alpha|=j}\frac{1}{\alpha
%!}\partial^\alpha H(0)y^\alpha$ is a homogeneous polynomial of
%degree $j,$ and $F_0(\imath)=G(\imath,\mu)+H(\imath).$

%We need some ``wave invariants'',  Kronecker's theorem.

%Why Gutzwiller trace formula does not work for non-degenerate
%critical point and resonances generated by closed orbits?

%Let us first put $E_0=0.$   In \cite{GuilleminUribe2007} it was
%assumed that potential has unique non-degenerate global minimum,
%$V(0)=0,$ at $x=0:$ $d=0$ in (\ref{V}). Moreover, it was assumed
%that,   for $\epsilon
%>0$ sufficiently small, $V^{-1}([0,\epsilon])$ is compact, then
%there is $h_0 >0$ such that for $h <h_0$ the spectrum of $P$ in a
%small interval, $[0,\delta],$ consists of a finite number of
%discrete eigenvalues, $E(h),$ and by Weyl's law
%\begin{equation}\label{Weyl}
%\sharp\,\{ E(h);\,\,0\leq E(h)\leq\delta\}=(2\pi h)^{-n}\left({\rm
%Vol}\,\{0\leq \xi^2 +V(x)\leq\delta\} +o(1)\right)
%\end{equation}

{\em Acknowledgements.} The author thanks the unknown referee for numerous comments and suggestions.

\section{Classical Birkhoff canonical form, proof of Lemma
\ref{lemma1}}\label{Section1}

Conjugating the Hamiltonian $p(x,\xi)=\xi^2+ V(x),$ with $V$ as in
(\ref{expansion}),
%$$V(x)=\sum_{j=1}^{n-d} u_j^2x_j^2 -\sum_{j=n-d+1}^n u_j^2x_j^2
%+{\mathcal O}(|x|^4)$$ (with only even terms present)
 by the linear
symplectomorphism
$$ x_i\mapsto u_i^{1/2}x_i,\,\,\xi_i\mapsto
u_i^{-1/2}\xi,\,\,i=1,\ldots, n,$$ one can assume without loss of
generality that $$ p=E_0+H_1+V_2\equiv E_0+\sum_{j=1}^{n-d}
u_j(\xi_j^2 +x_j^2) +\sum_{j=n-d+1}^n u_j(\xi_j^2-x_j^2)
+V_2(x_1^2,\ldots,x_n^2),$$ where $V_2(s_1,\ldots, s_n)={\mathcal O}
(|s|^2).$ We denote $H=H_1+V_2.$

Then resonances can essentially (see Introduction) be viewed as an
eigenvalue problem for $P$ after the complex scaling
\begin{equation}\label{tildes}
x_j=e^{i\pi/4} \tilde{x}_j,\,\,\tilde{x}_j\in\rz,\,\,n-d+1\leq j\leq
n.
\end{equation}
 The principal
symbol of the scaled operator becomes
$\tilde{p}(\ldots,\tilde{x},\ldots,\tilde{\xi})=E_0+\tilde{H}_1+\tilde{V}_2,$
where the new $\tilde{H}=\tilde{H}_1+\tilde{V}_2$ is equal to
$$H(x_1,\ldots, x_{n-d}, e^{i\pi /4}
\tilde{x}_{n-d+1},\ldots, e^{i\pi /4} \tilde{x_n}, \xi_1,\ldots,
\xi_{n-d}, e^{-i\pi /4} \tilde{\xi}_{n-d+1},\ldots, e^{-i\pi /4}
\tilde{\xi_n}).$$ With  $x_j=e^{i\pi/4}\tilde{x}_j,$
$\xi_j=e^{-i\pi/4}\tilde{\xi}_j,$
 for $n-d+1\leq j\leq n,$ we have
 $\xi_j^2-x_j^2=(\tilde{\xi}_j^2 +\tilde{x}_j^2)/i,$
 %
%---------------------------------------------------------------
%
%Here we need just formal scaling or restriction of $H(x,\xi)$
%\begin{equation}\label{scaledspace}
%(\rz^{n-d}\times e^{i\pi /4}\rz^{d})\times (\rz^{n-d}\times e^{-i\pi
%/4}\rz^{d}). \end{equation} Writing $x_j=e^{i\pi/4}\tilde{x}_j,$
%$\xi_j=e^{-i\pi/4}\tilde{\xi}_j,$
% for $n-d+1\leq j\leq n,$ we have
% $$\xi_j^2-x_j^2=e^{-i\pi/2}\tilde{\xi}_j^2 -
% e^{i\pi/2}\tilde{x}_j^2=(\tilde{\xi}_j^2 +\tilde{x}_j^2)/i.$$
%
%----------------------------------------------------------------
%
 and omitting the tildes we get
 $$  H(x,\xi)=\sum_{j=1}^{n-d} u_j(\xi_j^2 +x_j^2)
+\sum_{j=n-d+1}^n \frac{1}{i}u_j(\xi_j^2+x_j^2)
+V_2(x_1^2,\ldots,x_n^2)$$ with all $x_j,\xi_j$ real, and which can be
identified with the restriction of the old $H$ to the IR-manifold $\Lambda\in\cz^{2n}.$ Then
one can follow Guillemin-Uribe \cite{GuilleminUribe2007} keeping in
mind that for $n-d +1\ldots\leq j\leq n,$ $u_j$ are exchanged by
$u_j/i.$

%{\tt In Kaidi-Kerdelhue Birkhoff transform was applied before
%Bargman.  I need Birkhoff transform on the Bargman side.}

We have
\begin{equation}\label{H1}
H_1=\sum_{j=1}^{n-d} u_j(\xi_j^2+ x_j^2) +\sum_{j=n-d+1}^n
\frac{1}{i}u_j(\xi_j^2+x_j^2).
\end{equation}
As in  \cite{GuilleminUribe2007} we introduce complex coordinates,
$z_j=x_j +i\xi_j,$ with real $x_j,\xi_j.$  In these coordinates
$x_j^2+\xi_j^2=z_j\overline{z}_j=|z_j|^2.$
 The Hamiltonian vector field
$$\nu=\sum_j\frac{\partial H_1}{\partial\xi_j}\frac{\partial}{\partial
x_j}-\frac{\partial H_1}{\partial x_j}\frac{\partial}{\partial
\xi_j}$$ becomes the vector field
$$\frac{2}{i}\sum_{j=1}^{n-d}u_j\left(z_j\frac{\partial}{\partial
z_j} -\overline{z}\frac{\partial}{\partial\overline{z}_j}\right) -2
\sum_{j=n-d+1}^{n}u_j\left(z_j\frac{\partial}{\partial z_j}
-\overline{z}\frac{\partial}{\partial\overline{z}_j}\right).$$ Then
the proof of  \cite{GuilleminUribe2007},  where $u_j$ for $n-d+1\leq
j\leq n$ are substituted by $ u_j/i,$ can be applied and we get
inductively that for $N=1,2,\ldots$ there exists a \neigh{},
${\mathcal O},$ of $x=\xi=0,$ and a complex \canform{},
$\kappa:\,\,{\mathcal O}\mapsto\cz^{2n}$ such that
\begin{equation}\label{classicalBnf}
\kappa^*H=\sum_{j=1}^N H_j+R_{N+1}+R_{N+1}',
\end{equation}
 where
\begin{enumerate}
         \item[a)] The $H_j$ are homogeneous polynomials of degree
         $2j$ of the form $H_j=h_j(x_1^2 +\xi_1^2,\ldots,x_{n}^2
         +\xi_{n}^2),$ with $H_1$ given in (\ref{H1}).
         \item[b)] $R_N$ is homogeneous of degree $2N$ and of the
         form $R_N=W_N+R_N^\sharp,$ where $W_N$ consists of the
         terms homogeneous of degree $2N$ in the Taylor series of
         $V(x_1^2,\ldots,x_n^2)$ at $x=0,$ and $R_N^\sharp$ is an
         artifact of the previous inductive steps.
         \item[c)] $R_N'$ vanishes to order $2N+2$ at the origin and
         is of the form $R_N'=V-\sum_{k=2}^NV_k+S_N,$ where $S_N$ is
         another artifact of the inductive process. In addition,
         $R_N'$ is even.
       \end{enumerate}
Using this induction argument Guillemin and Uribe show that one can
read off from the $H_j$'s the first $N$ terms in the Taylor
expansion of $V(s_1,\ldots,s_n)$ at $s=0.$ This argument is
invariant under complex scaling. This achieves the proof of Lemma
\ref{lemma1}. \phantom{=}\hfill\qed\vspace{0.5cm}

 Recalling the tildes introduced by (\ref{tildes}) and letting $N$ tend
to infinity in (\ref{classicalBnf}) we obtain the classical Birkhoff
normal form
$$\sum_{j=1}^\infty \tilde{H}_j(x_1^2 +\xi_1^2,\ldots,x_{n-d}^2 +\xi_{n-d}^2,\tilde{x}_{n-d+1}^2
+\tilde{\xi}_{n-d+1}^2,\ldots, \tilde{x}_{n}^2 +\tilde{\xi}_{n}^2)$$
with $\tilde{H}_1$ as in (\ref{H1}). Then  after scaling back to
$\rz^n\times\rz^n$  we get the classical Birkhoff normal form as in
(\ref{F0}):
$$F_0=\sum_{j=1}^\infty H_j(\xi_1^2+x_1^2 ,\ldots,\xi_{n-d}^2+x_{n-d}^2 ,
\xi_{n-d+1}^2-x_{n-d+1}^2,\ldots, \xi_{n}^2-x_{n}^2 ),$$ with
$$H_1=\sum_{j=1}^{n-d} u_j(\xi_j^2 +x_j^2) +\sum_{j=n-d+1}^n
u_j(\xi_j^2-x_j^2).$$ The construction of the quantum Birkhoff
normal form (\ref{Bnf}) is well known (see for example
\cite{Sjostrand1992}).

%\bibliography{coulomb}

\end{document}